# A Fully-Distributed Asynchronous Approach for Multi-Area Coordinated Network-Constrained Unit Commitment

Yamin Wang, Lei Wu, Jie Li

*Electrical and Computer Engineering Department, Clarkson University, Potsdam, NY 13699 USA*

Email: yamwang@clarkson.edu, lwu@clarkson.edu, jieli@clarkson.edu

**Abstract**— This paper discusses a consensus-based alternating direction method of multipliers (ADMM) approach to solve the multi-area coordinated network-constrained unit commitment (NCUC) problem in a distributed manner. Due to political and technical difficulties, it is neither practical nor feasible to solve the multi-area coordination problem in a centralized fashion, which requires full access to all data of individual areas. In comparison, in the proposed fully-distributed approach, local NCUC problems of individual areas can be solved independently, and only limited information is exchanged among adjacent areas for facilitating multi-area coordination. Furthermore, as traditional ADMM can only guarantee convergence for convex problems, this paper discusses several strategies to mitigate oscillations, enhance convergence performance, and derive good-enough feasible solutions, including: (i) A tie-line power flow based area coordination strategy is designed to reduce the number of global consensus variables; (ii) Different penalty parameters $\rho$ are assigned to individual consensus variables and are updated via certain rules during the iterative procedure, which would reduce the impact of initial values of $\rho$ on convergence performance; (iii) Heuristic rules are adopted to fix certain unit commitment variables for avoiding oscillations during the iterative procedure; and (iv) An asynchronous distributed strategy is studied, which solves NCUC subproblems of small areas multiple times and exchanges information with adjacent areas more frequently within one complete run of slower NCUC subproblems of large areas. Numerical cases illustrate effectiveness of the proposed asynchronous fully-distributed NCUC approach, and investigate key factors that would affect its convergence performance.

# Nomenclature

**Indices:**

| | |
|---|---|
| $d, i$ | Indices of loads and units |
| $g$ | Index of global variables |
| $k, n/n'$ | Indices of iterations and areas |
| $t, \tau$ | Indices of hours |
| $m, j$ | Indices of buses |
| $mj$ | Index of a transmission line connecting bus $m$ and bus $j$ |

**Variables:**

| | |
|---|---|
| $F_i()$ | Fuel cost of unit $i$ |
| $I_{i,t}$ | Unit commitment variable of unit $i$ at time $t$ |
| $P_{i,t}$ | Real power output of unit $i$ at time $t$ |
| $PL_{mj,t}$ | Real power flow from bus $m$ to bus $j$ at time $t$ |
| $ON_{i,t}/OFF_{i,t}$ | ON/OFF time counter of unit $i$ at time $t$ |
| $SU_{i,t}/SD_{i,t}$ | Startup/shutdown cost of unit $i$ at time $t$ |
| $Z_{g,t}$ | Global consensus variable $g$ at time $t$ |
| $\theta_{m,t}$ | Voltage angle of bus $m$ at time $t$ |

**Constants:**

| | |
|---|---|
| $C_{su_i}/C_{sd_i}$ | Startup/shutdown cost coefficient of unit $i$ |
| $I_{i,0}$ | Initial commitment status of generator $i$ |
| $C_{nl_i}$ | No-load cost of unit $i$ |
| $NT$ | Number of hours under study |
| $ON_{i,0}, OFF_{i,0}$ | Initial ON/OFF time counter of unit $i$ |
| $P_{d,t}$ | Demand level of load $d$ at time $t$ |
| $PL_{mj}^{max}$ | Capacity limit of a line connecting buses $m$ and $j$ |
| $P_i^{min}, P_i^{max}$ | Minimum/maximum capacity of unit $i$ |
| $T_{on,i}/T_{off,i}$ | Minimum ON/OFF time limit of unit $i$ |
| $UR_i/DR_i$ | Ramp up/down rate limit of unit $i$ |
| $X_{mj}$ | Reactance of a line connecting buses $m$ and $j$ |

**Sets and Vectors:**

| | |
|---|---|
| $B_n, B_n(m)$ | Set of buses /set of buses connected to bus $m$ of area $n$ |
| $D_n(m)$ | Set of loads located at bus $m$ of area $n$ |
| $GL$ | Set of global variables |
| $I_n, P_n$ | Vector of unit commitment / real power output variables for generators located in area $n$ |
| $L_n$ | Set of tie-lines in area $n$ |
| $N, T$ | Set of areas /time periods |
| $U_n, U_n(m)$ | Set of units located in area $n$ /set of units located at bus $m$ in area $n$ |
| $\theta_n$ | Vector of voltage angles for buses located in area $n$ |
| $\lambda_n$ | Vector of Lagrangian multipliers corresponding to area $n$ |

# I. Introduction

Network-constrained unit commitment (NCUC) is one of the fundamental decision-making tools in power market operation, which minimizes total operation cost of generators while satisfying prevailing system and unit operation constraints. Traditionally, independent system operators (ISOs) and regional transmission organizations (RTOs) solve NCUC in a centralized manner, for optimizing asset utilization in their own service territory. However, in recognizing the intensified interconnection of regional electricity networks and the deeper penetration of renewable energy resources, coordinated operation of multiple geographical areas is highly desired to achieve the overall energy security and economic efficiency. Indeed, as indicated in [1], roughly 50% of the time power flows of tie-lines between NYISO and ISONE were in the wrong direction (i.e., electricity flows from a high price region to a low price region). Recently, several interregional interchange-scheduling projects have been conducted by ISOs and RTOs for improving the joint economic efficiency and achieving more cohesive boundary coordination, such as Inter-Regional Interchange Scheduling by ISONE and NYISO [1], Interchange Optimization by PJM and MISO [2], and Coordinated Transaction Scheduling by PJM and NYISO [3].

Ideally, a centralized solution to multi-area coordination, by gathering all data from individual areas, may achieve better overall energy security and economic efficiency of the entire interconnected electricity infrastructure. However, such a centralized solution framework is unlikely to be practical because of political and technical difficulties. In fact, individual ISOs and RTOs may not be willing to disclose their actual financial information, system topology, or control regulations. As a result, the multi-area coordination problem has been studied in a distributed manner [4]-[6]. [4] analyzed the challenges for defining decision-making strategies that are fair to individual system operators of an interconnected electricity infrastructure. A distributed procedure for scheduling multi-area optimal reactive power flows was discussed in [5]. Optimal power scheduling for an interconnected multi-area power system with cross-border trading was investigated in [6] while considering wind power uncertainty.

This paper focuses on solving the multi-area coordinated NCUC problem in a fully-distributed fashion without the central coordination. That is, in distributed NCUC (DNCUC), smaller-scale local NCUC problems of individual areas are solved in a distributed manner, and only limited information (i.e., power flows of tie-lines) is exchanged with adjacent areas during an iterative distributed optimization procedure.

Various approaches have been explored in literature to solve power system operation problems in a distributed fashion.

(1) Lagrangian Relaxation (LR): In this method, coupling constraints among different areas are relaxed via Lagrangian multipliers, which makes the derived problem separable [7]-[9]. Many further established distributed algorithms are based on LR.
(2) Auxiliary Problem Principle (APP): APP can be interpreted as solving a sequence of auxiliary problems involving augmented LR, in order to improve convergence performance of standard LR [10]. It was applied to distributed optimal power flow (OPF) problems by duplicating variables in coupling constraints [11]-[12]. In [6], an APP based multi-area UC with wind power uncertainty was discussed.
(3) Alternating Direction Method of Multipliers (ADMM): ADMM, theoretically well suited for convex problems [13]-[17], has also been extended to solve nonconvex operation problems of power systems in a distributed manner. A compact formulation of ADMM was discussed in [15] for solving AC-OPF problems. In [17], OPF was formulated as a semi-definite programming model and solved in a distributed fashion via ADMM.
(4) Analytical Target Cascading (ATC): An ATC-based decentralized UC approach was designed in [18] for the coordination between the ISO and distribution grid operators (DISCOs). Moreover, when coefficient used to update penalty parameters $\rho$ is set to 1, the ATC algorithm in [18] can be viewed as a standard ADMM. It was further applied in [19] for solving decentralized UC of large-scale power systems with a star topology.

However, directly applying above algorithms to solve multi-area coordinated NCUC problem could encounter convergence issues, as NCUC with binary unit commitment variables is highly non-convex and non-derivable. Indeed, it may converge to a local solution or even fail in finding a feasible solution to the original non-convex problem.

Existing distributed algorithms have been modified for solving mixed-integer programming problems in [20]-[22]. An ADMM based heuristic method for mixed-integer quadratic programming problems was discussed in [20], which sacrifices solution accuracy for achieving a reduced computational time. [21] discussed several heuristics and refinements, including fixing binary variables and adjusting values of penalty parameters, to mitigate oscillations of ADMM for solving distributed UC. [22] applied an ADMM based alternating optimization procedure (AOP) to solve the robust decentralized UC problem with wind power uncertainty, which is based on the facts that the original UC is naturally separable when voltage angles of boundary buses are fixed and further becomes convex when all binary variables are fixed. However, when fixing voltage angles of boundary buses to values obtained from the relaxed UC problem, local UC problems in the first iteration of AOP may be infeasible. Furthermore, the number of iterations between fix-theta and fix-binary procedures in AOP is always limited in practical cases, which indicates that AOP may fail to perform sufficient explorations for identifying good-enough solutions. In addition, [23] discussed several heuristic strategies to improve convergence performance of progressive hedging (PH) algorithm for stochastic mixed-integer problems. These strategies include selecting effective $\rho$ for independent decision variables, designing near-convergence detection rules, and fixing certain decision variables with cyclic behavior. [24] further evaluated solution quality of the PH algorithm by computing lower bound of optimal objective function. This paper solves the multi-area coordinated NCUC problem via a fully-distributed consensus-based ADMM approach. Different from the authors' previous work [29] which solves a convex DC-OPF problem while convergence property of many distributed algorithms could be exactly satisfied, distributed UC studied in this paper is nonconvex and directly applying traditional distributed algorithms could face significant convergence issues because of binary unit commitment variables. In fact, they may converge to a local solution or even fail in finding a feasible solution to the original non-convex non-differentiable UC problem.

Major contributions of the paper are fourfold.

(1) Different from the bus voltage angle based decomposition strategy in our previous work [29], a tie-line power flow based area decomposition strategy is designed for the consensus based ADMM algorithm. This strategy involves fewer global consensus variables which could improve convergence performance.

(2) Key factors that would aggravate oscillations when directly applying traditional ADMM on DNCUC are analyzed. Several strategies are discussed for accelerating convergence performance and achieving good-enough feasible solutions, including adopting distinct $\rho$ to individual consensus variables, updating $\rho$ iteratively to mitigate the impact of initial settings on $\rho$, and designing heuristic rules to fix oscillating binary variables.

(3) As computational time in each iteration of synchronous distributed algorithms is restricted by the slowest sub-problem, an asynchronous distributed procedure is studied.

(4) Numerical case studies compare the proposed strategies with APP [6], AOP [22], and PH [23] to illustrate its effectiveness.

The rest of this paper is organized as follows. Section II presents formulation of the ADMM based multi-area coordinated NCUC problem. Section III discusses several strategies to mitigate oscillations and enhance convergence performance of the distributed NCUC. Numerical case studies and conclusions are presented in Sections IV and V, respectively.

## II. ADMM based NCUC formulation

### A. NCUC Problem for a Single Area

The NCUC problem for a single area $n$ is formulated as a mixed-integer linear programming (MILP) problem (1), which optimizes hourly generation scheduling for supplying electricity loads [27]. Objective function (1a) is to minimize the total operation cost, including energy production cost $F(P_{i,t})$, no-load cost $C_{nl_i} \cdot I_{it}$, startup cost $SU_{i,t}$, and shutdown cost $SD_{it}$. Constraint (1b) represents generation capacity limits. Constraints (1c)-(1d) describe minimum ON/OFF time limits of generators, which represent that generator $i$ must be online/offline for at least $T_{on,i}/T_{off,i}$ hours before it can be switched off/on [39]. Constraints (1e)-(1f) enforce startup and shutdown costs. Constraints (1g)-(1h) ensure ramping up and down limits. Constraint (1i) represents nodal load balance equations which enforce the total power injection to a bus $m$ equals to the total power withdraw form that bus. Constraint (1j) calculates DC power flow on a transmission line connecting bus $m$ and bus $j$. Constraint (1k) is the DC transmission network security limit. Constraint (1l) enforces voltage phase angle of reference bus to be 0. An illustrative NCUC example for a simple 2-area 4-bus system is provided in Appendix.

$$\min\left(\sum_{i \in U_n} \sum_{t \in T}\left[F(P_{i,t}) + C_{nl_i} \cdot I_{i,t} + SU_{i,t} + SD_{i,t}\right]\right) \quad (1a)$$

$$s.t. P_i^{min} \cdot I_{i,t} \leq P_{i,t} \leq P_i^{max} \cdot I_{i,t}, \forall t \in T, \forall i \in U_n \quad (1b)$$

$$\sum_{t=1}^{UT_i}(1 - I_{i,t}) = 0, \quad \text{where } UT_i = \max\{0, \min[NT, (T_{on,i} - ON_{i,0}) \cdot I_{i,0}]\}, \forall i \in U_n$$

$$\sum_{\tau=t}^{t+T_{on,i}-1} I_{i,\tau} \geq T_{on,i} \cdot (I_{i,t} - I_{i,(t-1)}), \quad \forall t = UT_i + 1, \cdots, NT - T_{on,i} + 1, \forall i \in U_n$$

$$\sum_{\tau=t}^{NT}[I_{i,\tau} - (I_{i,t} - I_{i,(t-1)})] \geq 0, \quad \forall t = NT - T_{on,i} + 2, \cdots, NT, \forall i \in U_n \quad (1c)$$

$$\sum_{t=1}^{DT_i} I_{i,t} = 0, \text{where } DT_i = max\{0, min[NT, (T_{off,i} - OFF_{i,0})] \cdot (1 - I_{i,0})\}, \forall i \in U_n$$

$$\sum_{\tau=t}^{t+T_{off,i}-1}(1 - I_{i,\tau}) \geq T_{off,i} \cdot (I_{i,(t-1)} - I_{i,t}), \forall t = DT_i + 1, \cdots, NT - T_{off,i} + 1, \forall i \in U_n$$

$$\sum_{\tau=t}^{NT}[1 - I_{i,\tau} - (I_{i,(t-1)} - I_{i,t})] \geq 0, \quad \forall t = NT - T_{off,i} + 2, \cdots, NT, \forall i \in U_n \quad (1d)$$

$$SU_{i,t} \geq C_{su_i} \cdot (I_{i,t} - I_{i,(t-1)}), \quad \forall t \in T, \forall i \in U_n \quad (1e)$$

$$SD_{i,t} \geq C_{sd_i} \cdot (I_{i,(t-1)} - I_{i,t}), \quad \forall t \in T, \forall i \in U_n \quad (1f)$$

$$P_{i,t} - P_{i,(t-1)} \leq UR_i \cdot I_{i,(t-1)} + P_i^{min} \cdot (I_{i,t} - I_{i,(t-1)}) + P_i^{max} \cdot (1 - I_{i,t}), \forall t \in T, \forall i \in U_n (1g)$$

$$P_{i,(t-1)} - P_{i,t} \leq DR_i \cdot I_{i,t} + P_i^{min} \cdot (I_{i,(t-1)} - I_{i,t}) + P_i^{max} \cdot (1 - I_{i,(t-1)}), \forall t \in T, \forall i \in U_n \quad (1h)$$

$$\sum_{i \in U_n(m)} P_{i,t} = \sum_{j \in B_n(m)} PL_{mj,t} + \sum_{d \in D_n(m)} P_{d,t}, \quad \forall m \in B_n, \forall t \in T \quad (1i)$$

$$PL_{mj,t} = \frac{1}{X_{mj}} (\theta_{m,t} - \theta_{j,t}), \quad \forall mj \in L_n, \forall t \in T \quad (1j)$$

$$PL_{mj}^{max} \leq PL_{mj,t} \leq PL_{mj}^{max}, \quad \forall mj \in L_n, \forall t \in T \quad (1k)$$

$$\theta_{ref,t} = 0, \quad \forall t \in T \quad (1l)$$

For the sake of discussion, the NCUC problem (1) for area $n$ is presented in a compact form as in (2).

$$\min \left( \sum_{i \in U_n} \sum_{t \in T} [F(P_{i,t}) + C_{nl_i} \cdot I_{i,t} + SU_{i,t} + SD_{i,t}] \right) \quad (2a)$$

$$s.t. \; \boldsymbol{g}_n(\boldsymbol{P}_n, \boldsymbol{\theta}_n, \boldsymbol{PL}_n, \boldsymbol{I}_n) \leq 0 \quad (2b)$$

$$\boldsymbol{h}_n(\boldsymbol{P}_n, \boldsymbol{\theta}_n, \boldsymbol{PL}_n, \boldsymbol{I}_n) = 0 \quad (2c)$$

Where $\boldsymbol{P}_n$, $\boldsymbol{\theta}_n$, and $\boldsymbol{PL}_n$ are respectively vectors of continuous variables for real power output, bus angle, and transmission line power flow associated with area $n$; $\boldsymbol{I}_n$ is set of binary unit commitment variables associated with area $n$; $\boldsymbol{g}_n(.)$ represents inequality constraints (1b)-(1h) and (1k) associated with area $n$; $\boldsymbol{h}_n(.)$ includes equality constraints (1i), (1j), and (1l) associated with area $n$.

## B. Distributed Multi-Area Coordinated NCUC Problem

The multi-area coordinated NCUC problem can be mathematically represented as in (3), which minimizes the total operation cost over all areas subject to operational and security constraints of individual areas as well as their coupling relationships.

$$\min \sum_{n \in N} \sum_{i \in U_n} \sum_{t \in T} [F(P_{i,t}) + C_{nl_i} \cdot I_{i,t} + SU_{i,t} + SD_{i,t}] \quad (3a)$$

$$s.t. \quad (2.b)-(2.c) \text{ for each area } n \quad (3b)$$

The multi-area coordinated NCUC model (3) presents a peculiar structure with respect to local variables and constraints. Specifically, objective function (3a) can be naturally split into objectives of individual areas, which are solely based on information of local generators. Constraints (1b)-(1h) in (2b) are also based on local variables of individual areas. On the other hand, areas are coupled with each other via power flows on tie-lines, which are naturally inseparable. For instance, for two boundary buses $j$ and $m$ of two areas as shown in Fig. 1, $PL_{jm,t}$ in constraints (1i) and (1k) indicates power flow exchange between areas $n$ and $n'$. In fact, DC power flows on tie-lines are calculated by voltage angles of boundary buses (1l). That is, $PL_{jm,t} = (\theta_{j,t} - \theta_{m,t})/X_{jm}$, where $\theta_{j,t}$ and $\theta_{m,t}$ are voltage angles of boundary buses $j$ and $m$ in areas $n$ and $n'$ at time $t$, respectively.

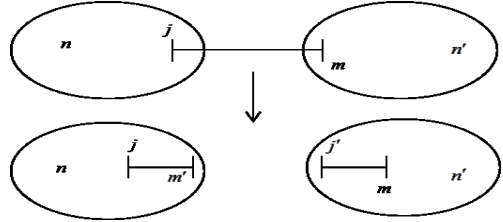

Fig. 1. Illustration of bus voltage angle based area decomposition strategy

To bring constraints (1i) and (1k) into a suitable form for distributed optimization, a widely used area decoupling approach is to duplicate voltage phase angles of boundary buses in one area to its adjacent areas [18], [29]-[32]. That is, as illustrated in Fig. 1, pseudo bus $m'$ in area $n$ is a duplication of bus $m$ in area $n'$, which can be viewed as a local copy of boundary bus $m$ in area $n$. Similarly, pseudo bus $j'$ in area $n'$ is a duplication of bus $j$ in area $n$. Same variables duplicated in different areas are forced to be equal via additional constraints $\theta_{m,t} = \theta_{m',t}$ and $\theta_{j,t} = \theta_{j',t}$. In equations (4), two global consensus variables $Z_{1,t}$ and $Z_{2,t}$ are introduced to link same bus angle variables duplicated in different areas.

$$\theta_{j,t} = Z_{1,t}, \theta_{m',t} = Z_{2,t}, \quad \forall t \in T \quad (4a)$$

$$\theta_{j',t} = Z_{1,t}, \theta_{m,t} = Z_{2,t}, \quad \forall t \in T \quad (4b)$$

Our previous work [29] indicated that a useful guidance for designing an efficient distributed algorithm is to avoid a large number of global variables. That is, a smaller number of global variables usually indicates fewer coupling variables, less information exchange, and a faster convergence rate. In this paper, we propose a tie-line power flow based area decoupling method by introducing duplicated tie-line power flows variables, instead of boundary bus voltage angles, into problem (3), in order to reduce the number of global consensus variables. That is, in the proposed tie-line power flow based area decoupling method, $\widetilde{PL}_{jm',t}$ and $\widetilde{PL}_{m'j,t}$ represent power flows of a same tie-line that are calculated separately by the two areas as in (5a) and (5b). The advantage of the proposed method is that only one global consensus variable $Z_t$ is needed to link the same tie-line power flow represented in different areas as in (5c), comparing with two global variables in (4a)-(4b).

$$\widetilde{PL}_{jm',t} = (\theta_{j,t} - \theta_{m',t})/X_{jm}, \forall t \in T \quad (5a)$$

$$\widetilde{PL}_{mj',t} = (\theta_{m,t} - \theta_{j',t})/X_{jm}, \forall t \in T \quad (5b)$$

$$\widetilde{PL}_{jm',t} = Z_t, \widetilde{PL}_{mj',t} = -Z_t, \quad \forall t \in T \quad (5c)$$

With above stated area decoupling method, the compact form of the multi-area coordinated NCUC problem (3) can be re-formulated as in (6), where $\overline{PL}_n$ represents power flow variables of local transmission lines in area $n$, $\widetilde{PL}_n$ indicates tie-line power flow variables in area $n$, and $\widetilde{PL}_{n,t}(w)$ represents the $w$th tie-line power flow variable in area $n$ at time $t$. Specifically, (6a) represents summation of objectives for individual areas, where $C_n(.)$ represents the objective function (1a) of area $n$. (6b) and (6c) respectively represent all local inequality and equality constraints in area $n$. (6d) represents that the same tie-line power flow calculated in different areas should be equal to a unique global consensus variable $Z$, where $G_t(n,w)=g$ represents the mapping from variable $\widetilde{PL}_{n,t}(w)$ onto the global consensus variable $Z_{g,t}$ at time $t$ (Equation (6d) is further illustrated via a 2-area 4-bus system in Appendix).

$$\min \sum_{n \in N} C_n(\boldsymbol{P}_n, \boldsymbol{\theta}_n, \boldsymbol{I}_n, \overline{\boldsymbol{PL}}_n, \widetilde{\boldsymbol{PL}}_n) \quad (6a)$$

$$s.t. \boldsymbol{g}_n(\boldsymbol{P}_n, \boldsymbol{\theta}_n, \boldsymbol{I}_n, \overline{\boldsymbol{PL}}_n, \widetilde{\boldsymbol{PL}}_n) \leq 0, \quad \forall n \in N \quad (6b)$$

$$\boldsymbol{h}_n(\boldsymbol{P}_n, \boldsymbol{\theta}_n, \boldsymbol{I}_n, \overline{\boldsymbol{PL}}_n, \widetilde{\boldsymbol{PL}}_n) = 0, \quad \forall n \in N \quad (6c)$$

$$\widetilde{PL}_{n,t}(w) = Z_{g,t}, \forall t \in T, \forall G_t(n,w) = g, \forall g \in \boldsymbol{GL} \quad (6d)$$

The multi-area coordinated NCUC formulation (6) is suitable for the consensus based ADMM procedure [28]-[30]. (7a) is the augmented Lagrangian function of (6a) after relaxing (6d) into the objective function, where $x_n$ represents all decision variables $\boldsymbol{P}_n, \boldsymbol{\theta}_n, \boldsymbol{I}_n, \boldsymbol{PL}_n, \widetilde{\boldsymbol{PL}}_n$ in area $n$, $x_n \in \chi_n$ indicates $x_n$ satisfying (6b)-(6c) in area $n$, $\lambda_{n,t}$ is set of Lagrangian multipliers corresponding to area $n$ at time $t$, $\lambda_{n,t}(w)$ is the $w$th variable of $\lambda_{n,t}$, and $\rho$ is a predefined positive augmented Lagrangian parameter. The $k$th iteration of consensus based ADMM procedure includes three steps of (7b)-(7d). (7b) for individual areas can be solved in parallel. (7c) means that a global variable $Z_{g,t}$ equals to average of all $\widetilde{PL}_{n,t}(w)$ that correspond to $Z_{g,t}$. (7d) represents the procedure of updating dual variable $\lambda_{n,t}(w)$. In the method of multipliers, a dual variable can be updated via a step size equal to the augmented Lagrangian parameter $\rho$. The algorithm terminates when primal and dual residuals (7e) in each area $n$ are sufficiently small [28]. The stopping criterion (7e) has been widely used and tested in literature [15], [29], [30] and [38], and has been proved in [28] to be effective for determining the convergence of ADMM in optimizing power system operation problems. Applicability of other stopping criteria in the proposed ADMM based DNCUC, such as those discussed in [46], will be analyzed in authors' future work.

$$\min_{x,z} L_\rho(x, z, \lambda) = C_n(x_n) + \sum_{\forall t \in T} \sum_{\forall w \in \lambda_{n,t}} \lambda_{n,t}(w) \cdot (\widetilde{PL}_{n,t}(w) - Z_{g,t})$$
$$+ \sum_{\forall t \in T} \sum_{\forall w \in \lambda_{n,t}} (\rho/2) \cdot \left\| \widetilde{PL}_{n,t}(w) - Z_{g,t} \right\|_2^2 \quad (7a)$$

$$x_n^{k+1} = \operatorname{argmin}_{x_n \in \chi_n} \begin{pmatrix} C_n(x_n) + \sum_{t \in T} \sum_{w \in \lambda_{n,t}} \lambda_{n,t}^k(w) \cdot \left( \widetilde{PL}_{n,t}(w) \right) + \\ \sum_{t \in T} \sum_{w \in \lambda_{n,t}} (\rho/2) \cdot \left\| \widetilde{PL}_{n,t}(w) - Z_{g,t}^k \right\|_2^2 \end{pmatrix}, \forall n \in N \quad (7b)$$

$$Z_{g,t}^{k+1} = \sum_{G_t(n,w)=g} \widetilde{PL}_{n,t}^{k+1}(w) / \sum_{G_t(n,w)=g} 1, \qquad \forall g \in GL, \forall t \in T \quad (7c)$$

$$\lambda_{n,t}^{k+1}(w) = \lambda_{n,t}^k(w) + \rho \cdot \left( \widetilde{PL}_{n,t}(w) - Z_{g,t} \right), \quad \forall w \in \lambda_{n,t}, \forall n \in N, \forall t \in T \quad (7d)$$

$$\|r_n^{k+1}\|_2^2 = \|\lambda_n^{k+1} - \lambda_n^k\|_2^2 \leq \varepsilon_1, \forall\ n \in N;\ \|s^{k+1}\|_2^2 = \rho \cdot \|Z^{k+1} - Z^k\|_2^2 \leq \varepsilon_2 \quad (7e)$$

The distributed optimization procedure for the multi-area coordinated NCUC problem is summarized as Algorithm 1.

**Algorithm 1:** ADMM based DNCUC

1. Initialize **Z** and $\lambda_n$ for all areas $n \in N$, and set iteration index $k$=1.
2. Execute one iteration of ADMM. That is, the following calculations for individual areas $n \in N$ are executed in parallel:
   - 2.1 Each area *n* solves its local NCUC problem (7b).
   - 2.2 Each area *n* exchanges its local solutions of tie-line power flows with adjacent areas.
   - 2.3 Each area *n* updates $Z_{gt}$ and $\lambda_{nt}(w)$ via (7c) and (7d), respectively.
3. If stopping criteria (7e) are satisfied for all areas, the algorithm terminates; Otherwise, set $k$=$k$+1 and go back to Step 2.

Proper settings on initial values of **Z** and $\lambda_n$ in Step 1 could affect convergence performance of Algorithm 1. In this paper, these initial values are determined via a continuous relaxation version of Algorithm 1. That is, we first set all **Z** and $\lambda_n$ as zeros, relax all binary variables in (7b) as continuous variables between 0 and 1, and execute Algorithm 1 to optimality. Final solutions of **Z** and $\lambda_n$ will be used as initial values in Step 1 of the ADMM based NCUC. In addition, as number of global variables **Z** plays a significant role on convergence performance of the consensus based ADMM algorithm [29], performance of the proposed tie-line power flow based area decoupling method and the voltage phase angle based area decoupling method will be compared via case studies in Section IV.

## III. Enhancement on the ADMM Based Distributed Multi-Area Coordinated NCUC

Since global convergence of ADMM can only be guaranteed for convex problems, directly applying Algorithm 1 on the multi-area coordinated NCUC may fail in converging to a feasible solution. This section focuses on strategies to enhance convergence performance of Algorithm 1, including modifications to parameter $\rho$, heuristic rules to fix certain unit commitment variables, and an asynchronous distributed strategy.

### A. Modifications to Parameter $\rho$

When directly applied to nonconvex NCUC problems, Algorithm 1 could fail to converge because (7e) may never be satisfied, as shown in Fig. 2(a). Numerical case studies show that failure of convergence is mainly caused by oscillation of binary unit commitment variables. As shown in Fig. 2(b), unit commitment solutions are oscillating among four states during the iterative procedure, where those four points in each column represent four different unit commitment combinations of three generators. One factor that may aggravate oscillations of binary variables is value of penalty parameter $\rho$. In fact, in one of our tests, the multi-area coordinated NCUC problem converges to a feasible solution with $\rho$=4 but oscillates with $\rho$=2.

Furthermore, as analyzed in [29], $\rho$ is also an important factor that influences convergence speed of Algorithm 1. Thus, a proper value of $\rho$ is important for mitigating oscillations and improving convergence performance.

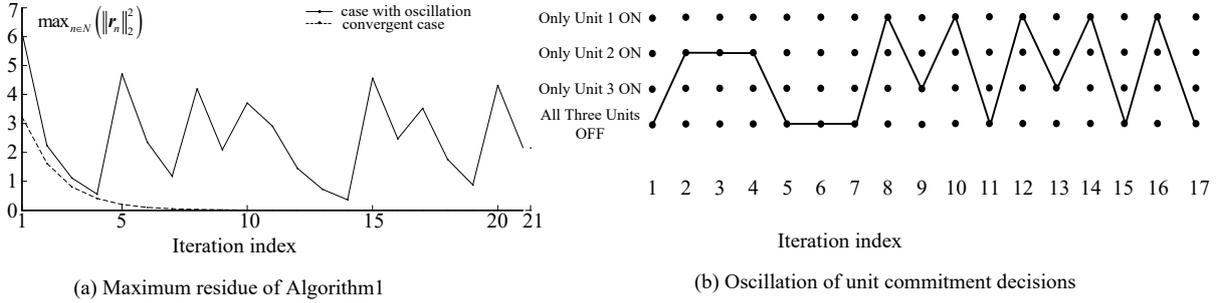

(a) Maximum residue of Algorithm1        (b) Oscillation of unit commitment decisions

Fig. 2. An illustrative example of oscillation

The authors in [21] regulate convergence performance of ADMM by heuristically adjusting value of a unique $\rho$ to update all dual variables. Specifically, [21] indicated that increasing value of $\rho$ could force settlement on a certain integer feasible solution, while decreasing $\rho$ could encourage more explorations on new integer feasible solutions. However, our tests show that solely adjusting $\rho$ may fail to perform enough explorations. In fact, although reducing value of $\rho$ may help the algorithm jump out of current local integer feasible solution, it is more likely to fall into a previously identified solution instead of exploring new ones.

Different from [21], in this paper, in order to improve convergence performance and make the algorithm less dependent on initial setting of $\rho$, two modifications are studied. (i) In observing that a unique parameter $\rho$ may not be the best choice for updating all $\lambda$ in (7d), each global variables $Z_{g,t}$ is given a distinct parameter $\rho_{g,t}$. (ii) Since parameter $\rho_{g,t}$ is usually set by experience and may not be universal for all practical cases, we update $\rho_{g,t}$ in each iteration via (8), where $\alpha^{incr}$, $\alpha^{decr}$, and $\mu$ are parameters larger than 1. The idea behind this penalty parameter updating scheme (8) is to keep one residual within a factor $\mu$ of the other, when they are both approaching to zero. That is, $\rho_{g,t}$ is increased by a factor of $\alpha^{incr}$ when primal residual $r_{g,t}$ is larger than $\mu$ times of dual residual $s_{g,t}$, which would accelerate the decrease of primal residual $r_{g,t}$ in next iteration; On the other hand, $\rho_{g,t}$ is decreased by a factor of $\alpha^{decr}$ when dual residual $s_{g,t}$ is much larger than primal residual $r_{g,t}$, in order to reduce dual residual $s_{g,t}$ faster in next iteration. Adjusting $\rho_{g,t}$ along the iterative procedure could make the algorithm less dependent on initial choice and provide proper values in each iteration to accelerate convergence. Indeed, as indicated in [28], although it is difficult to prove convergence of ADMM when $\rho$ varies along iterations, (8) often works well in many practical cases. Numerical tests in Section IV indicate that proposed modifications to parameter $\rho$ can accelerate convergence speed and make performance less dependent on initial choice of $\rho$. In addition, the convergence theory for ADMM with fixed-$\rho$ could be still applicable to ADMM with (8), if $\rho$ becomes fixed after a finite number of ADMM iterations.

$$\rho_{g,t}^{k+1} = \begin{cases} \alpha^{incr}\rho_{g,t}^{k}, & if \ r_{g,t}^{k} > \mu \cdot s_{g,t}^{k} \\ \rho_{g,t}^{k}/\alpha^{decr}, & if s_{g,t}^{k} > \mu \cdot r_{g,t}^{k} \\ \rho_{g,t}^{k}, & otherwise \end{cases} \quad (8)$$

## B. Heuristic Rules to Fix Certain Unit Commitment Variables

Binary fixing strategies have been discussed in literature [21] to accelerate convergence. That is, if value of a binary variable remains unchanged for a certain number of iterations, this binary variable will be fixed to this value. However, fixing too many binary variables may lead to premature convergence, i.e., ADMM could easily converge to a sub-optimal solution with low solution quality. In turn, in order to obtain solutions with acceptable accuracy, the algorithm in [21] needs to explore multiple feasible solutions before termination, which significantly increase computational complexity. In addition, based on our extensive observations, failure of convergence is mainly caused by oscillation of unit commitment variables in the iterative procedure, while fixing variables whose solutions remain unchanged [21] may not contribute to mitigating convergence issue. Thus, different from [21], this paper proposes to fix binary variables that present oscillation behavior. In this paper, a unit commitment variable is identified as an oscillating variable if solution to this variable changes frequently in a number of consecutive ADMM iterations.

Three spatial and temporal heuristic rules are studied in this section to determine priority of generators in individual areas of different time periods whose unit commitment variables are to be fixed. This would help mitigate oscillation issue and facilitate convergence performance.

(1) The first rule is to identify the sequence of hours in which oscillation occurs and the priority of binary variables that need to be fixed. In a certain ADMM iteration when oscillation in a consecutive multi-hour period is identified, binary variables of the first hour in that oscillating period will have the highest priority to be fixed. For instance, if oscillations occur at hours 5-7 and 14-16, hour 5 and hour 14 have the highest priority and their unit commitment variables will be fixed. This rule is based on our observation that oscillations occurred in consecutive hours are mainly caused by minimum ON/OFF time constraints of generating units, while once oscillation in the first hour is mitigated by fixing those binary variables, oscillation in the entire consecutive multi-hour period might be effectively resolved.

(2) The second rule is that an area with the most tie-line connections will have the highest priority to fix binary variables. This rule is set based on the observation that, oscillation of unit commitment variables in one area will result in the vibration of tie-line power flows and further influence unit commitment variables in its adjacent areas. Once unit commitment variables in this area are fixed, related oscillations in adjacent areas could also be mitigated. In turn, fixing binary variables for the area with the largest number of tie-lines first would help mitigate oscillations in as many adjacent areas as possible.

(3) The last rule is to define priority list of generators (*pg*) that unit commitment variables need to be fixed, which is calculated by the full-load average production cost (9) [33]. Specifically, expensive generators have larger *pg*, and in turn will be given a lower priority to be fixed. This is consistent with power system operation practice that cheaper units are usually served as base units with relative stable unit commitment schedules. On the other hand, unit commitment schedules of expensive units are more fluctuating and their impacts on the objective function are more significant. In turn, it would be better to fix unit commitment variables of cheaper units and leave commitment variables of expensive units to be optimized in later iterations.

$$pg_i = F(P_i^{max})/P_i^{max} \tag{9}$$

Finally, Algorithm1 is modified as Algorithm 2 with above enhanced strategies. Specifically, in Step 4 of Algorithm 2, the oscillation detection procedure begins after 30 iterations of ADMM algorithm, and a unit commitment variable is identified as an oscillating variable if its value changes more than 3 times in 10 consecutive iterations. In addition, one ADMM iteration in this paper refers to a complete run of all subproblems for all areas, i.e., Steps 2.1-2.3 of Algorithm 1 or Steps 2.1-2.4 of Algorithm 2. Furthermore, binary variables that have been fixed in previous ADMM iterations will remain fixed and not be released in later iterations. Practical studies show that number of oscillating generators in each area is usually limited and new sets of oscillating time periods during the iterative solution procedure are also limited. Thus, the algorithm is expected to converge to a solution after fixing several groups of binary variables. In addition, [23] also proposed to fix certain complicating variables for mitigating oscillations in the PH algorithm when solving stochastic UC problems, while the idea of fixing certain key binary variables to mitigate oscillations of distributed algorithms for non-convex optimization problems is of a general interest. It is also worth mentioning that the heuristic fixing strategy may not be able to entirely solve oscillation issue in the distributed NCUC problem, which deserves a strict mathematical proof to verify in future work.

**Algorithm 2:** Enhanced ADMM based DNCUC
1. Initialize **Z** and $\lambda_n$ for each area $n \in N$, and set iteration index $k$=1, $k_{fix}$=1;
2. Execute one iteration of ADMM. That is, following calculations for individual areas $n \in N$ are executed in parallel:
    2.1 Each area $n$ solves its local NCUC problem (7b).
    2.2 Each area $n$ exchanges its local solutions of tie-line power flows with adjacent areas.
    2.3 Each area $n$ updates $Z_{gt}$ and $\lambda_{nt}(w)$ via (7c) and (7d), respectively.
    2.4 Each area $n \in N$ updates $\rho_g$ via (8).
3. If stopping criteria (7e) for all areas are satisfied, the algorithm terminates; Otherwise, go to Step 4.
4. If $k$ is larger than 30 and ($k$ - $k_{fix}$)>10, the following procedure is executed to detect oscillating variables; Otherwise, we go to Step 5.
    4.1 Detect time periods in which oscillation occurs, and record the set of hours with the highest priority based on the first rule, denoted as $T^{pri}$
    4.2 **For** each $t$ in $T^{pri}$
    4.3    Based on the second rule, determine the set of areas with the highest priority, denoted as $N_t^{pri}$.
    4.4    **For** each $n$ in $N_t^{pri}$
    4.5       Based on the third rule, determine the set of generators in $N_{n,t}^{pri}$ with the highest priority, denoted as $G_{g,n,t}^{pri}$.
    4.6    **END**
    4.7 **END**
    4.8 Fix all identified oscillating unit commitment variables $G_{g,n,t}^{pri}$ to 1, and set $k_{fix}$=$k$.
5. Set $k$=$k$+1 and go back to Step 2.

## C. Asynchronous Fully-Distributed Algorithm

As computational time in each iteration of the synchronous distributed algorithm is limited by the slowest area, an asynchronous DNCUC is further studied, in which small areas optimize their local NCUC sub-problems and exchange tie-line power flow solutions with their adjacent areas multiple times rather than waiting for large areas. That is, (7c) and (7d) can be performed immediately after certain related areas finish updating (7b).

In literature, [35] proposed a randomized asynchronous ADMM algorithm in which different components of the network are allowed to wake up randomly and perform local updates, while the rest of the network stands still. In the randomized asynchronous ADMM algorithm [35], no coordinator or global clock is needed. [36] further indicated that convergence rate of asynchronous ADMM for convex problems is $O(1/k)$. [25] designed a asynchronous distributed ADMM for consensus optimization, in which partial barrier $S$ and bounded delay $\tau/T$ are used to control asynchrony. [37] proposed a distributed ADMM procedure in which individual agents operate asynchronously based on the coloring scheme of network. [40] designed an asynchronous algorithm assembling a multi-block ADMM. Instead of relying on probabilistic control or restrictive assumptions, the problem instance in [40] only makes standard convex-analytic regularity assumptions. In addition, [26] and [34] indicated that asynchronous ADMM could converge to the set of Karush-Kuhn-Tucker condition points for possibly non-convex problems in the form of (6a)-(6d), if parameters of the algorithm are properly chosen.

In this paper, the asynchronous ADMM framework follows the scheme proposed in [25]. The detailed asynchronous procedure and the convergence analysis for several typical problems, including network average consensus problem, can be referred to [25]. In asynchronous ADMM, partial barrier $S$ represents that at least $S$ areas need to update (7b) before (7c)-(7d) can be performed, and bounded delay $\tau/T$ means that (7b) for each area has to be executed and updated with adjacent areas at least $\tau$ times in every $T$ ADMM iterations. Fig. 3 illustrates the synchronous distributed ADMM procedure as well as the asynchronous distributed ADMM procedure with $S=1$ and $\tau/T=1/2$ for an interconnected three-area system. As shown in Fig. 3(b), faster areas such as Area 2 can update their local solutions more frequently in asynchronous ADMM.

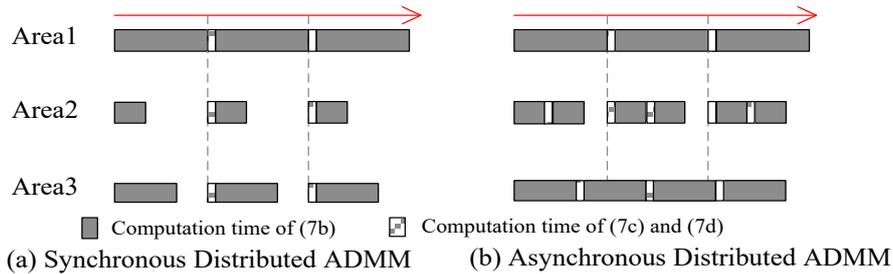

(a) Synchronous Distributed ADMM  (b) Asynchronous Distributed ADMM

Fig. 3 Illustrative example of synchronous and asynchronous distributed ADMM procedures

In addition, in the asynchronous case, the binary variable fixing strategy is performed when a global ADMM iteration terminates, i.e., after the slowest area completes one local NCUC calculation, and values to be fixed are based on the latest NCUC results of individual areas.

## IV. Numerical Case Studies

Four sets of tests are discussed in this section. The first test analyzes cases in which the traditional ADMM could converge to a local solution, while performance of the proposed tie-line power flow based area decomposition strategies is compared with the APP based decentralized algorithm [6]. The second test discusses cases in which directly applying ADMM fails to converge to a local solution, while convergence performance of the proposed heuristic methods is compared with the AOP algorithm [22]. The third test further compares the heuristic strategies proposed in this paper with the methods in [23]-[24] via a large-scale 4-area 472-bus system. The fourth test illustrates the asynchronous DNCUC via a complicated 8-area system. In all tests, local NCUC problems (7b) are implemented in MATLAB with YALMIP and solved via Gurobi.

### A. Tests on the Proposed Enhancements to the Traditional ADMM

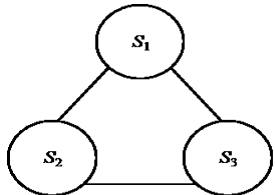

Fig. 4. A 3-area 72-bus system

In the 3-area system shown in Fig. 4, each area is a modified IEEE 24-bus system with 10 generators, 17 loads, and 34 transmission lines. Hourly load profile for the entire system is shown in Table I. $\varepsilon_1$ and $\varepsilon_2$ are set to $10^{-3}$, and relative MIP gap for solving (7b) is set as 0.1%. Two Cases $A_1$ and $A_2$ are discussed. Case $A_1$ compares computational performance of bus voltage angle and tie-line power flow based decomposition strategies discussed in Section II.B. In Case $A_2$, the ADMM based DNCUC is compared with the APP based decentralized algorithm [6]. No oscillation occurs in these two cases.

Table I System Load Profile (% of the Total Installed Capacity)

| Hour | 1 | 2 | 3 | 4 | 5 | 6 | 7 | 8 |
|---|---|---|---|---|---|---|---|---|
| Demand | 54.4 | 54.8 | 55.5 | 54.1 | 54.5 | 56.2 | 60.8 | 62.2 |
| Hour | 9 | 10 | 11 | 12 | 13 | 14 | 15 | 16 |
| Demand | 65.4 | 72.4 | 80.1 | 82.9 | 85.1 | 85.4 | 87.1 | 89.6 |
| Hour | 17 | 18 | 19 | 20 | 21 | 22 | 23 | 24 |
| Demand | 89.9 | 86.5 | 86.1 | 83.3 | 83.3 | 81.5 | 68.5 | 68.5 |

***Case $A_1$***: In this case, the proposed tie-line power flow based area decomposition strategy (i.e., Strategy 2) is compared with the bus voltage angle based one (i.e., Strategy 1). Algorithm 1 is applied to compare the two strategies. Results are shown in Table II.

Table II shows that Strategy 2 converges faster than Strategy 1 for various values of $\rho$. The reason is that the number of global variables in Strategy 2 is 3×24(h)=72, which is smaller than 6×24(h)=144 in Strategy 1. In addition, in this case study, Strategy 1 with different $\rho$ values converges to a same local optimal solution with the total operation cost of $3,041,006, and Strategy 2 with different $\rho$ values also converges to a same local optimal solution with the total operation cost of $3,026,568. In summary, the proposed tie-line power flow based decomposition Strategy 2, by reducing 50% of global variables as compared to Strategy 1, presents a better convergence performance than Strategy 1.

Table II Number of Iterations for the Two Decomposition Strategies

|  | $\rho=2$ | $\rho=4$ | $\rho=8$ | $\rho=16$ | $\rho=32$ |
|---|---|---|---|---|---|
| Strategy1 | 56 | 54 | 58 | 62 | 123 |
| Strategy2 | 23 | 21 | 22 | 23 | 52 |

***Case $A_2$***: [6] and [11] illustrated that APP exhibits a better convergence performance than ADMM for non-convex problems. In this case study, the APP based DNCUC algorithm is also applied to the system in Fig. 4 and compared with the proposed Algorithm 2. Furthermore, in an ideal situation, a central controller can solve the multi-area coordinated NCUC problem with all data from individual areas, and its optimal solution is used as a benchmark. Specifically, difference in total operation costs between the decentralized algorithm and the ideal centralized method (10) is used as a metric to evaluate solution quality of the proposed fully-distributed algorithm. In (10), $Tcost_{dec}$ and $Tcost_{cen}$ are total operation costs of decentralized and centralized methods, respectively.

$$Diff=|Tcost_{dec} - Tcost_{cen}|/Tcost_{cen} \times 100\% \tag{10}$$

The APP based algorithm solves the multi-area coordinated NCUC problem in 7.3s after 19 iterations, with the total operation cost of $3,026,536. In comparison, as shown in Table III, the best computational performance of Algorithm 2 is achieved with $\rho_{g,t}=4$, which converges in 7.4s after 18 iterations with the total operation cost of $3,026,568. Furthermore, total operation cost of the centralized algorithm is $3,019,131, which is obtained in 3s. Thus, *Diff* of both Algorithm 2 and APP based algorithm are around 0.25%. In summary, both algorithms derive high quality solutions while presenting similar computational performance.

Table III Number of Iterations with Modifications to $\rho$ in Algorithm 2

| Initial $\rho_{g,t}$ | 2 | 4 | 8 | 16 | 32 | Centralized |
|---|---|---|---|---|---|---|
| Iteration # | 20 | 18 | 19 | 22 | 28 | N/A |
| Total computational time (s) | 8.2 | 7.4 | 7.9 | 9.3 | 11.9 | 3.0 |

## B. Tests on Proposed Heuristic Methods for Fixing Binary Variables

As discussed in Section III, directly applying Algorithm 1 to the multi-area coordinated NCUC may fail in converging to a feasible solution. In this section, the proposed heuristic methods for fixing binary variables are further tested on the multi-area system shown in Fig. 5. The test system consists of 4 areas, while each area is a modified IEEE 24-bus system. Directly applying Algorithm 1 on this system fails to converge to a solution. This system is further solved by Algorithm 2 with heuristic rules to fix certain oscillating binary variables. Another heuristic method, AOP algorithm proposed in [22], is also applied for comparison. The relative MIP gap for solving (7b) is set as 0.1% in this test.

In this case, oscillations happen at hours 8-10 in areas $S_2$, $S_3$, and $S_4$ during the ADMM procedure. Fig.2 (b) shows oscillations of unit commitment solutions of three generators in area 2 at hour 8. Based on the first rule in Section III.B, unit commitment variables at hour 8 have the highest priority to be fixed. In addition, based on the second rule in Section III.B, oscillating generators 1-3 in area $S_2$, which has more tie-lines than areas $S_3$ and $S_4$, have the highest priority to be fixed. Furthermore, generator 1 in area $S_2$ has the lowest $pg$ among generators 1-3. Thus, based on the third rule in Section III.B, unit commitment variable of generator 1 in area $S_2$ at hour 8 is first fixed as 1. After this single fixing action is applied, Algorithm 2 converges to a solution after 20 iterations.

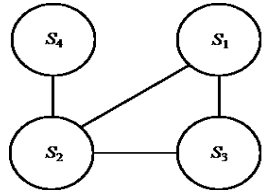

Fig. 5. A 4-area 96-bus system

In comparison, if we do not follow the first rule and randomly choose unit commitment variables to be fixed, for instance, unit commitment variables at hour 9 in area $S_2$, Algorithm 2 will not converge until an additional group of unit commitment variables at hour 8 is fixed. If we do not follow the second rule and randomly choose area $S_4$ as the first area to be dealt with, oscillations in areas $S_2$ and $S_3$ still persist even oscillation in area $S_4$ is diminished. These additional case studies show that the proposed heuristic rules could provide valuable guidance on determining priorities of unit commitment variables to be fixed and in turn improve convergence performance.

In this case, ON/OFF schedules of all generators from Algorithm 2 and the centralized approach are compared in Table IV. Results in Table IV indicate that at certain hours, number of committed units in Algorithm 2 is larger than that of the centralized solution. For instance, at hour 1, 18 generators are committed in Algorithm 2 while only 16 are ON in the centralized solution. This phenomenon can be explained as follows. During the iterative procedure of Algorithm 2, after solving (7b) area $S_1$ may require X MW from area $S_2$ to avoid turning on one additional generator in area $S_1$ with $P_{min}$ much larger than X MW. In a same iteration, area $S_2$ may be facing the same situation. That is, solution to (7b) for area $S_2$ may require Y MW from area $S_1$ to avoid turning on one additional generator in area $S_2$ with $P_{min}$ much larger than Y MW. If these two areas cannot reach an agreement during the consensus procedure (7b)-(7d), oscillation will happen; otherwise, if they could reach an agreement, our observations show that the final solution of Algorithm 2 would most likely turn on generators in both areas, although turning on only one unit might be the global optimal solution for the entire interconnected system. Thus, in the decentralized algorithm, more generators might be turned ON. In order to avoid this phenomenon and reduce the number of unnecessarily committed generators, in the proposed fixing strategy, we only deal with unit commitment variables that perform oscillations, which could reduce the chance of over-commitment issues as many unit commitment variables do not present oscillation behavior along the ADMM iterative procedure. Furthermore, the proposed heuristic strategy only fixes a very limited number of oscillating unit commitment variables that have the highest priority as described in Section III.B, which could further reduce the chance of over-commitment issues. For instance, the second heuristic rule designed in Section III.B is set to first fix unit commitment variables of the area with the largest number of tie-lines, which could also help mitigate oscillations of its adjacent areas and in turn reduce the number of unit commitment variables in those areas that need to be fixed in later iterations. In fact, *Diff* of 0.17% for Algorithm 2 indicates that solution obtained by Algorithm 2 is very close to the centralized one.

In addition, *Diff* of the AOP algorithm is 0.72% which is larger than that of Algorithm 2. The total number of ADMM iterations is 56 for the AOP algorithm, while is 35 for Algorithm 2. In each iteration of the AOP algorithm, (7b) is solved with all binary variable fixed and in turn computing time for solving (7b) in AOP would be much shorter than that of Algorithm 2. In fact, the total computational time of Algorithm 2 is 22.4s, while is 12.8s for the AOP algorithm. Thus, compared to the AOP algorithm, Algorithm 2 shows advantage of higher solution accuracy at the cost of longer computational time.

Table IV Solution Obtained by Algorithm 2 and the Centralized Algorithm

|   | Units | Algorithm 2 (24 hours) | Centralized Algorithm (24 hours) |
|---|---|---|---|
| $S_1$ | 1 | 000000011111111111111111 | 000000001111111111111111 |
|   | 2 | 000000000111111111111100 | 000000000111111111111100 |
|   | 3-5,7-8 | 111111111111111111111111 | 111111111111111111111111 |
|   | 6 | **11111**111111111111111111 | **00000**111111111111111111 |
|   | 9 | 000000000000**000000000**000 | 000000000000**111111111**000 |
|   | 10 | 000000000000000000000000 | 000000000000000000000000 |
| $S_2$ | 1 | 000000011111111111111111 | 000000011111111111111111 |
|   | 2 | 000000000**1111111111**1100 | 000000000**0000000000**0000 |
|   | 3 | 000000000**1111111111**1100 | 000000000**0000000000**0000 |
|   | 4 | 000000000000011100000000 | 000000000000000000000000 |
|   | 5-7,9-10 | 111111111111111111111111 | 111111111111111111111111 |
|   | 8 | **111111**111111111111111111 | **000000**111111111111111111 |
| $S_3$ | 1 | 000000011111111111111111 | 000000011111111111111111 |
|   | 2 | 000000000111111111111100 | 000000000111111111111100 |
|   | 3-8 | 111111111111111111111111 | 111111111111111111111111 |
|   | 9 | 000000000000000**000**0000000 | 000000000000000**111**0000000 |
|   | 10 | 000000000000000000000000 | 000000000000000000000000 |
| $S_4$ | 1 | 000000011111111111111111 | 000000011111111111111111 |
|   | 2 | 000000000111111111111100 | 000000000111111111111100 |
|   | 3-8 | 111111111111111111111111 | 111111111111111111111111 |
|   | 9-10 | 000000000000000000000000 | 000000000000000000000000 |

## C. Comparison with heuristic strategies for the PH algorithm [23]-[24]

In this section, we discuss differences and similarities between heuristic strategies for PH algorithm in [23]-[24] and the proposed methods via the following two sub-sections. All cases in this section are tested on a large-scale 4-area 472-bus system. The interconnection topology of these four areas is the same as that in Fig. 5, while each area is replaced by a modified IEEE 118-bus system [29]. For the sake of discussion, strategies proposed in [23] are represented via variables defined in this paper. For instance, global consensus variable $\widetilde{PL}_{n,t}(w)$ is used to represent complicating variables mentioned in [23].

### *C1: Methods for Accelerating the Convergence Speed*

[23] proposed three methods to accelerate convergence performance of the PH algorithm. In this subsection, these three methods are first described briefly, followed by discussions on their potential applicability to the proposed distributed ADMM algorithm. The proposed

strategies in Section III.A are further tested on a 4-area 472-bus system. In this section, thresholds of stopping criteria $\varepsilon_1$ and $\varepsilon_2$ in (7e) are set as 0.05; relative MIP gap for solving (7b) is set as 0.1%; $\alpha^{inc}$, $\alpha^{decr}$, and $\mu$ in (8) are respectively set as 2, 2, and 5.

(1)  The first method proposed in [23] is to select a proper $\rho$ for each decision variable. It is suggested in [23] that an effective $\rho$ value of a variable should be close in magnitude to its cost coefficient in objective function. Based on this rule, distinct values of $\rho$ are assigned to different decision variables in [23]. However, this method may not be directly applicable to the proposed distributed ADMM algorithm. The reason is that in the distributed NCUC problem (7a)-(7e), $\rho$ corresponds to each tie-line power flow variable $\widetilde{PL}_{n,t}(w)$, which is not explicitly included in the objective function and does not possess a cost coefficient.

Alternatively, in this paper we propose to assign a distinct parameter $\rho_{g,t}$ for each global variable $Z_{g,t}$ and update value of $\rho_{g,t}$ in each iteration according to (8), although their initial values could be the same. These two modifications can make the algorithm less dependent on initial choice of $\rho_{g,t}$ and provide proper values along the iterative process to accelerate convergence. In fact, $\rho_{g,t}$ modification procedures proposed in [23] and this paper are similar in terms that they both intend to design variable-dependent $\rho_{g,t}$ values, by recognizing that a unified $\rho$ value for all variables may not be a proper choice for updating all $\lambda$ in (7d).

Effectiveness of the proposed modifications to $\rho$ in Section III.A is tested via a large scale 4-area 472-bus system, and results are shown in Table V. Specifically, when an improper initial setting of $\rho_{g,t}$ =20 is chosen, number of iterations of Algorithm 1 is 482 which is much larger than 102 of Algorithm 2. In addition, numbers of iterations of Algorithm 2 with respect to different initial values of $\rho_{g,t}$ are not significantly different. Results in Table V again verify that the proposed modifications to $\rho$ can make the algorithm less dependent on initial value of $\rho$ and also accelerate convergence speed. Results in Tables III and V also show that, computational time of the distributed algorithm is longer than that of the centralized solution. The main reason is that the distributed algorithm could need tens or even hundreds of iterations before it can converge to a feasible solution. Thus, although the complexity for solving (7a) in each area has been significantly reduced compared to the original centralized problem, the total computing time of the distributed algorithm with multiple iterations could be longer than that of the centralized solution.

Table V Computational performance for the 4-area 472-bus system with the proposed modifications to $\rho$

| Initial $\rho_{g,t}$ | 20 | 40 | 100 | 150 | 200 | Centralized |
|---|---|---|---|---|---|---|
| # of Iterations of Algorithm 1 | 482 | 420 | 243 | 184 | 216 | N/A |
| Total computational time (s) | 5591 | 4893 | 2843 | 1755 | 2529 | 390 |
| # of Iterations of Algorithm 2 | 102 | 98 | 82 | 73 | 112 | N/A |
| Total computational time (s) | 1198 | 1178 | 984 | 872 | 1340 | 390 |

(2)  The second method proposed in [23] is the design of near-convergence detection rules. [23] indicated that the PH algorithm empirically yields significant violation reductions in the early iterations to meet stopping criteria, while remaining iterations play a fine-tuning role to further drive already small violations to zero. Thus, it might be possible to terminate the iterative procedure earlier without losing significant solution accuracy, which could consequently improve computational efficiency. Indeed, [23] designed a simple near-convergence detection rule (11) based on primal complicating variables, which is used to terminate the algorithm earlier with an acceptable solution accuracy.

$$ts_{g,t} = \frac{\left(c_{max} \cdot max_{\forall G_t(n,w)=g} \widetilde{PL}_{n,t}(w)\right)}{\left(c_{min} \cdot min_{\forall G_t(n,w)=g} \widetilde{PL}_{n,t}(w)\right)} \leq \varepsilon_3, \forall g \in \boldsymbol{GL}, \forall t \in \boldsymbol{T} \tag{11}$$

where $c_{max}$ and $c_{min}$ are cost coefficients of corresponding complicating variables in the objective function.

ADMM presents similar convergence characteristics as PH. Indeed, it may only take ADMM a few tens of iterations to derive results with modest accuracy that are acceptable for practical engineering applications. To help illustrate this phenomenon, relationship between $max_{n \in N}(\|r_n\|_2^2)$ and the number of iterations for the 4-are 472-bus system is shown in Fig. 6. As illustrated in Fig. 6, after about 65 iterations, $max_{n \in N}(\|r_n\|_2^2)$ is reduced from 25.9 to 0.05, while it needs another 75 iterations to further reduce max($r$) from 0.05 to 0.01.

The near-convergence detection rule (11) defined in [23] cannot be directly applied in the proposed distributed ADMM based NCUC approach, because global tie-line power flow variables $\widetilde{PL}_{n,t}(w)$ are not included in the objective function and as a result we do not have $c_{max}$ and $c_{min}$ information for these global variables. Indeed, if $c_{max}$ equals to $c_{min}$, (11) basically measures discrepancy among duplicated complicating variables, which is similar to primal residual $r$ defined in (7e). Thus, an effective way in ADMM is to assign relatively large values of $\varepsilon_1$ and $\varepsilon_2$ in (7e) for terminating the algorithm earlier while not losing significant solution accuracy, which could be analogous to (11) proposed in [23]. By leveraging computational efficiency and solution accuracy, based on our extensive tests, proper values of $\varepsilon_1$ and $\varepsilon_2$ are in the range of 0.05 to 0.07 in distributed NCUC studies.

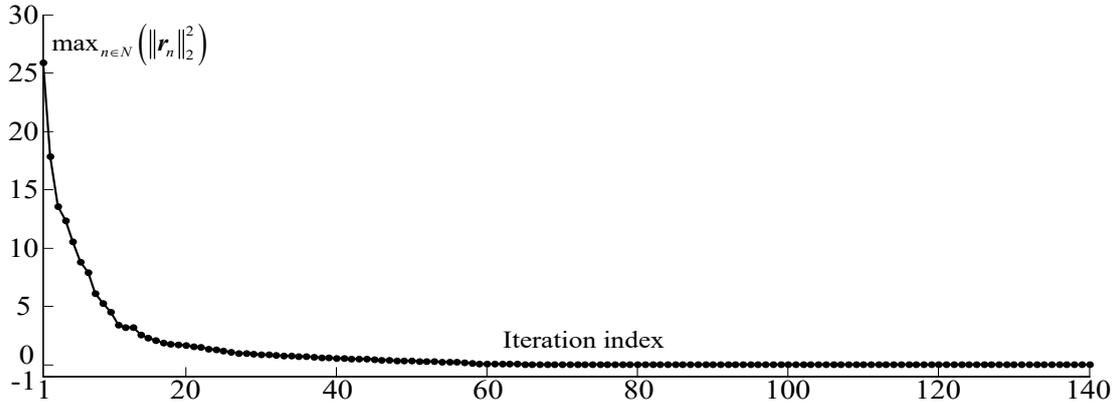

Fig.6 maximum value of $r$ versus iteration numbers for the 4-are 472-bus system

(3) [23] indicated that certain $\widetilde{PL}_{n,t}(w)$ variables may converge to a consensus $Z_{g,t}$ in only a few iterations while others may need more iterations to achieve convergence. Thus, the third method proposed in [23] for accelerating convergence is to fix $\widetilde{PL}_{n,t}(w)$ once they have converged to a consensus $Z_{g,t}$. This method is applied to the distributed ADMM algorithm, and test results on the 4-area 472-bus system are shown in Table VI.

Table VI Comparison of Algorithm 2 with/without Fixing Global Variables

| Initial $\rho_{g,t}$ | 20 | 40 | 100 | 150 | 200 |
|---|---|---|---|---|---|
| # of Iterations while fixing $\widetilde{PL}_{n,t}(w)$ | 101 | 98 | 82 | 72 | 112 |
| # of Iterations without fixing $\widetilde{PL}_{n,t}(w)$ | 102 | 98 | 82 | 73 | 112 |

As shown in Table VI, with different initial values of $\rho_{g,t}$, fixing $\widetilde{PL}_{n,t}(w)$ does not significantly impact number of iterations for the ADMM based DNCUC algorithm. Indeed, for $\rho_{g,t}$ =40, 100, and 200, numbers of ADMM iterations of the two methods are exactly the same. This phenomenon can be explained as follows. In this case study, when certain $\widetilde{PL}_{n,t}(w)$ variables converge to a solution, further iterations for driving other variables to a consensus will not significantly change values of global variables that have converged already. That is, values of global variables $\widetilde{PL}_{n,t}(w)$ that have already converged are stable in further iterations even though they are not explicitly fixed. For instance, in Fig. 7(a), $\widetilde{PL}_{2,1}(1)$ is fixed to the converged value after 40 iterations while in Fig. 7(b), $\widetilde{PL}_{2,1}(1)$ is left unfixed. As shown in Fig. 7(b), driving $\widetilde{PL}_{2,19}(1)$ to a consensus in further iterations does not significantly change the value of unfixed $\widetilde{PL}_{2,1}(1)$. It is noted that the $\widetilde{PL}_{n,t}(w)$ fixing procedure could be useful for cases in which further iterations for slow-convergent variables will result in significant re-update of already converged variables. Such situations will be further investigated in our future work.

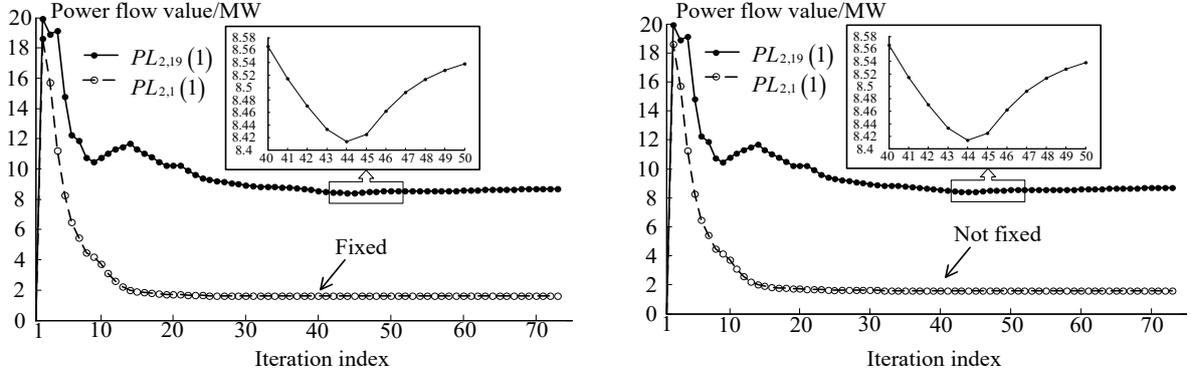

(a) Tie-line power flow vs. iteration number with $PL_{n,t}(w)$ fixed    (b) Tie-line power flow vs. iteration number without $PL_{n,t}(w)$ fixed

Fig. 7 Illustration of Algorithm 2 with/without fixing global variables

*C2: Methods for Mitigating the Risk of Non-convergence*

The method proposed in [23] for mitigating oscillations is to detect cyclic behavior of $\lambda_{n,t}(w)$ in (7d). That is, values of $\lambda_{n,t}(w)$ in multiple iterations are compared, and once the cyclic behavior of $\lambda_{n,t}(w)$ is detected, complicating variable $\widetilde{PL}_{n,t}(w)$ is fixed as the maximum value that appears in those cyclic iterations.

In comparison, in the NCUC problem, failure of convergence is mainly caused by oscillation of unit commitment variables. Thus, in this paper, we detect cyclic behavior of binary unit commitment variables instead of $\lambda_{n,t}(w)$, and adopt heuristic rules discussed in Section III.B to fix unit commitment variables of prioritized units that present cyclic behavior.

Effect of the heuristic rules for mitigating oscillations is further tested on the 4-area 472-bus system. As the 4-area 472-bus system studied in Section *C1* does not present cyclic behavior, in this subsection, load profile is modified to trigger a cyclic oscillation behavior. Oscillations in this case happen at hours 9-11 and hours 18-19 in areas $S_2$, $S_3$, and $S_4$ during the ADMM procedure. Based on the first rule discussed in Section III.B, unit commitment variables at hours 9 and 18 have the highest priority to be fixed. In addition, unit commitment variables of generators 1-2 in each of the three areas $S_2$, $S_3$, and $S_4$ perform oscillations at hour 9, and unit commitment variables of generators 3-5 in each of the three areas $S_2$, $S_3$, and $S_4$ perform oscillations at hour 18. Thus, based on the second rule discussed in Section III.B, oscillating generators in area $S_2$, which has more tie-lines than areas $S_3$ and $S_4$, have the highest priority to be fixed at both hours 9 and 18. Furthermore, generator 1 in area $S_2$ has a lower *pg* than generator 2, and generator 3 in area $S_2$ has the lowest *pg* among generators 3-5. Thus, based on the third rule discussed in Section III.B, unit commitment variables of generator 1 in area $S_2$ at hour 9 and generator 3 in area $S_2$ at hour 18 are first fixed to 1. As shown in Fig. 8, after applying this variable fixing step, Algorithm 2 converges to a solution after 17 ADMM iterations. This case study again verifies that the proposed heuristic rules can efficiently mitigate oscillation in the ADMM procedure.

In summary, both [23] and this paper focus on improving convergence performance and mitigating oscillations of decomposition based solution approaches for nonconvex optimization problems with binary variables. As discussed above, there are certain similarities between the methods discussed in [23] and the strategies proposed in this paper. However, some heuristic methods discussed in [23] may not be directly applicable to the ADMM based

DNCUC problem discussed in this paper. Additional studies on the 4-area 472-bus system further show that the heuristic strategies proposed in this paper are effective for enhancing convergence performance of the ADMM based distributed NCUC problem.

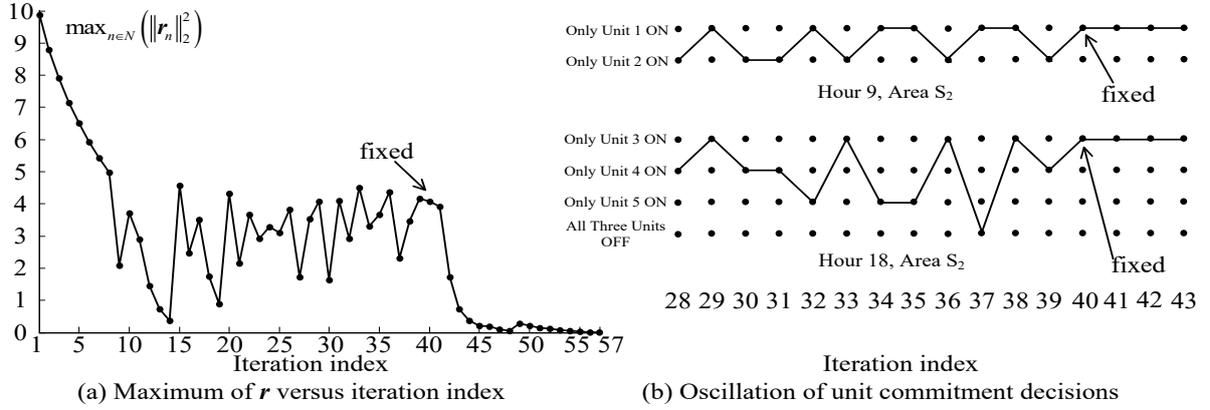

(a) Maximum of *r* versus iteration index    (b) Oscillation of unit commitment decisions

Fig. 8 Illustration on the effectiveness of the proposed fixing procedure

### D. Asynchronous Decentralized Algorithm

In this section, an 8-area system in Fig. 9 is studied. Area $S_1$ is a modified IEEE 118-bus system, and each of areas $S_2$-$S_8$ is a modified IEEE 24-bus system. Computing time for solving NCUC of $S_1$ is around 15s with the relative MIP gap of 0.5%, while is about 1s for solving each NCUC of $S_2$-$S_8$. Algorithm 2 is applied in this case study, in which $\rho_{g,t}$ is set as 8, while $\alpha^{inc}$, $\alpha^{decr}$, and $\mu$ in (8) are respectively set as 2, 2, and 5.

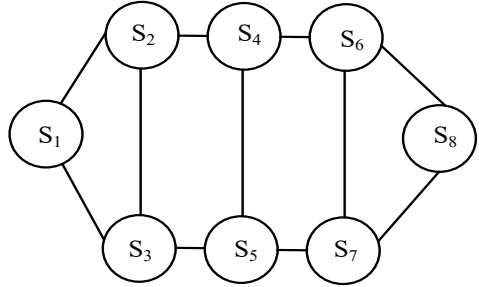

Fig. 9. An 8-area system

In Algorithm 1 and Algorithm 2, each Z in (7c) corresponds to one tie-line connecting two adjacent areas. Thus, in the synchronous procedure, partial barrier *S* equals to 2 which means (7c) and (7d) will be performed when both adjacent areas finish their local optimization (7b). In comparison, in the asynchronous ADMM *S* is set as 1, which means (7c) and (7d) could be updated if one area has finished the local optimization (7b). Bounded delay $\tau/T$ in this paper represents that for two adjacent areas, the slower area has to use information provided by the faster area at least $\tau$ times in every $T$ ADMM iterations.

The impact of $\tau/T$ on the convergence performance is tested via different values including 2/3, 1/2, 1/3, 1/4, and 1/8. Results are shown in Table VII. Numbers of iterations in Table VII are total numbers of iterations of the slowest area $S_1$, and the total computational time is based on the real wall clock time. As shown in Table VII, numbers of iterations and total computational times with $\tau/T$=2/3 and 1/2 are smaller than those of the synchronous procedure, which indicates an improvement in computational efficiency. Smaller $\tau/T$ means more self-updates of faster systems (i.e., $S_2$-$S_8$ in this case). Specifically, when $\tau/T$=1/2, $S_2$-$S_8$ calculate (7b) and update solutions twice while waiting for $S_1$ to finish its own optimization task (7b). This procedure could accelerate convergence of (7e). However, more updates of faster areas do not necessarily mean a better convergence performance because smaller $\tau/T$ also indicates that information from slower areas is incorporated into faster areas less frequently. As shown in Table VII, when $\tau/T$ takes values of 1/3 and 1/4, numbers of iterations and total computational time are all larger than the synchronous method.

Table VII Computational Performance of the Asynchronous Algorithm

| $\tau/T$ | Syn | 2/3 | 1/2 | 1/3 | 1/4 | Centralized |
|---|---|---|---|---|---|---|
| Number of Iterations | 54 | 52 | 49 | 60 | 82 | N/A |
| Total Computational Time (s) | 821.9 | 789.4 | 746.3 | 912.6 | 1253.8 | 253 |

In summary, results in Table VII show that, for non-convex NCUC problems, asynchronous ADMM could converge to a local optimal solution with acceptable accuracy. In addition, convergence performance of asynchronous ADMM depends on bounded delay $\tau/T$, while total computational time of asynchronous ADMM could be reduced when a proper bounded delay $\tau/T$ is set, i.e., computational time is reduced by 9.2% with $\tau/T$ of 1/2. Convergence issue of asynchronous ADMM for non-convex problems as well as the strategies for further improving its computational benefits will be analyzed in our future works.

## V. Conclusion

This paper discusses ADMM-based fully-decentralized optimization algorithms for solving multi-area coordinated NCUC problems. A tie-line power flow based area decomposition structure is designed for the consensus based ADMM algorithm. Our case studies indicate that this decomposition strategy presents a better convergence performance over traditional bus voltage angle based decomposition strategy. Since convergence of ADMM is only guaranteed for convex problems, directly applying ADMM on the multi-area coordinated NCUC may fail to converge and oscillations of binary variables could occur. Modifications to parameter $\rho$ and heuristic rules to fix oscillating binary variables are discussed in this paper. Numerical case studies indicate that they could significantly improve convergence performance while achieving good-enough solutions as compared to other methods in literature.

Asynchronous ADMM is also applied on the multi-area coordinated NCUC problem, in which faster small-scale areas update their optimal solutions and exchange necessary information with their adjacent areas multiple times before slower large-scale areas finish their local optimization tasks. Results show that the asynchronous algorithm could help improve computational efficiency, but too frequent updates in small-scale systems do not necessarily mean better convergence performance.

As a tighter UC formulation can potentially improve computational efficiency and solution accuracy [41]-[45], it could also contribute to an improved performance for solving distributed NCUC problems, although it may not directly resolve inherent convergence issue of ADMM based distributed NCUC approaches as studied in this paper. The impact of tighter UC formulations on computational performance of distributed NCUC will be investigated in authors' future work. Furthermore, conditions under which certain variables may be un-fixed and method to adjust values of variables [47]-[48] are also worth investigating. In addition, future work will also analyze convergence issue of asynchronous ADMM for non-convex problems, strategies for further improving computational benefits of asynchronous algorithms, and the over-commitment issue of distributed algorithms.

## Appendix: Illustration of proposed distributed NCUC via a 2-area 4-bus system

In this Appendix, an illustrative 2-area 4-bus system shown in Fig. 10(a) is used to describe the detailed formulation of NCUC problem (1) and the tie-line power flow based decomposition strategy proposed in Section II.B.

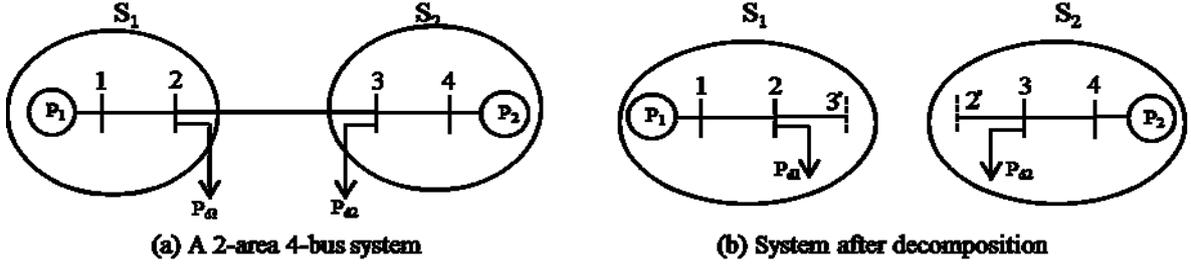

(a) A 2-area 4-bus system  (b) System after decomposition

Fig. 10. A 2-area 4-bus system

The centralized formulation of the NCUC problem for the 2-area 4-bus system in Fig. 10(a) is presented as in (A).

$$\min \sum_{t=1}^{NT}\{[F_1(P_{1,t}) + C_{nl1} \cdot I_{1,t} + SU_{1,t} + SD_{1,t}] + [F_2(P_{2,t}) + C_{nl2} \cdot I_{2,t} + SU_{2,t} + SD_{2,t}]\} \quad (A1)$$

$$s.t. \quad P_1^{min} \cdot I_{1,t} \leq P_{1,t} \leq P_1^{max} \cdot I_{1,t}, \quad \forall t \in T \quad (A2)$$

$$\sum_{t=1}^{UT_1}(1 - I_{1,t}) = 0, \text{ where } UT_1 = max\{0, min[NT, (T_{on,1} - ON_{1,0}) \cdot I_{1,0}]\}$$

$$\sum_{\tau=t}^{t+T_{on,1}-1} I_{1,\tau} \geq T_{on,1} \cdot (I_{1,t} - I_{1,(t-1)}), \quad \forall t = UT_1 + 1, \cdots, NT - T_{on,1} + 1$$

$$\sum_{\tau=t}^{NT}[I_{1,\tau} - (I_{1,t} - I_{1,(t-1)})] \geq 0, \quad \forall t = NT - T_{on,1} + 2, \cdots, NT \quad (A3)$$

$$\sum_{t=1}^{DT_1} I_{1,t} = 0, \quad \text{where } DT_1 = max\{0, min[NT, (T_{off,1} - OFF_{1,0})] \cdot (1 - I_{1,0})\}$$

$$\sum_{\tau=t}^{t+T_{off,1}-1}(1 - I_{1,\tau}) \geq T_{off,1} \cdot (I_{1,(t-1)} - I_{1,t}), \forall t = DT_1 + 1, \cdots, NT - T_{off,1} + 1$$

$$\sum_{\tau=t}^{NT}[1 - I_{1,\tau} - (I_{1,(t-1)} - I_{1,t})] \geq 0, \quad \forall t = NT - T_{off,1} + 2, \cdots, NT \quad (A4)$$

$$SU_{1,t} \geq C_{su_1} \cdot (I_{1,t} - I_{1,(t-1)}); \quad SD_{1,t} \geq C_{sd_1} \cdot (I_{1,(t-1)} - I_{1,t}) \quad \forall t \in T \quad (A5)$$

$$P_{1,t} - P_{1,(t-1)} \leq UR_1 \cdot I_{1,(t-1)} + P_1^{min} \cdot (I_{1,t} - I_{1,(t-1)}) + P_1^{max} \cdot (1 - I_{1,t}), \forall t \in T \quad (A6)$$

$$P_{1,(t-1)} - P_{1,t} \leq DR_1 \cdot I_{1,t} + P_1^{min} \cdot (I_{1,(t-1)} - I_{1,t}) + P_1^{max} \cdot (1 - I_{1,(t-1)}), \forall t \in T \quad (A7)$$

$$P_{1,t} - PL_{12,t} = 0; \quad \theta_{1,t} = 0, \quad \forall t \in T \quad (A8)$$

$$PL_{12,t} = (\theta_{1,t} - \theta_{2,t})/X_{12}; \quad PL_{12}^{min} \leq PL_{12,t} \leq PL_{12}^{max}; \quad \forall t \in T \quad (A9)$$

$$PL_{12,t} - PL_{23,t} = P_{d1,t}, \quad \forall t \in T \quad (A10)$$

$$P_2^{min} \cdot I_{2,t} \leq P_{2,t} \leq P_2^{max} \cdot I_{2,t}; \quad \forall t \in T \quad (A11)$$

$$\sum_{t=1}^{UT_2}(1 - I_{2,t}) = 0, \quad \text{where } UT_2 = max\{0, min[NT, (T_{on,2} - ON_{2,0}) \cdot I_{2,0}]\}$$

$$\sum_{\tau=t}^{t+T_{on,2}-1} I_{2,\tau} \geq T_{on,2} \cdot (I_{2,t} - I_{2,(t-1)}), \quad \forall t = UT_2 + 1, \cdots, NT - T_{on,2} + 1$$

$$\sum_{\tau=t}^{NT}[I_{2,\tau} - (I_{2,t} - I_{2,(t-1)})] \geq 0, \quad \forall t = NT - T_{on,2} + 2, \cdots, NT \quad (A12)$$

$$\sum_{t=1}^{DT_2} I_{2,t} = 0, \quad \text{where } DT_2 = max\{0, min[NT, (T_{off,2} - OFF_{2,0})] \cdot (1 - I_{2,0})\}$$

$$\sum_{\tau=t}^{t+T_{off,2}-1}(1 - I_{2,\tau}) \geq T_{off,2} \cdot (I_{2,(t-1)} - I_{2,t}), \forall t = DT_2 + 1, \cdots, NT - T_{off,2} + 1$$

$$\sum_{\tau=t}^{NT}[1 - I_{2,\tau} - (I_{2,(t-1)} - I_{2,t})] \geq 0, \quad \forall t = NT - T_{off,2} + 2, \cdots, NT \quad (A13)$$

$$SU_{2,t} \geq C_{su_2} \cdot (I_{2,t} - I_{2,(t-1)}); \quad SD_{2,t} \geq C_{sd_2} \cdot (I_{2,(t-1)} - I_{2,t}), \quad \forall t \in T (A14)$$

$$P_{2,t} - P_{2,(t-1)} \leq UR_2 \cdot I_{2,(t-1)} + P_2^{min} \cdot (I_{2,t} - I_{2,(t-1)}) + P_2^{max} \cdot (1 - I_{2,t}), \forall t \in T (A15)$$

$$P_{2,(t-1)} - P_{2,t} \leq DR_2 \cdot I_{2,t} + P_2^{min} \cdot (I_{2,(t-1)} - I_{2,t}) + P_2^{max} \cdot (1 - I_{2,(t-1)}), \forall t \in T (A16)$$

$$P_{2,t} + PL_{34,t} = 0; PL_{34,t} = (\theta_{3,t} - \theta_{4,t})/X_{34}; \quad PL_{34}^{min} \leq PL_{34,t} \leq PL_{34}^{max}, \forall t \in T (A17)$$

$$PL_{23,t} - PL_{34,t} = P_{d2,t}; PL_{23,t} = (\theta_{2,t} - \theta_{3,t})/X_{23}; PL_{23}^{min} \leq PL_{23,t} \leq PL_{23}^{max}, \forall t \in T (A18)$$

Objective function (A1) is to minimize total operation cost of the entire system, which is composed of two areas $S_1$ and $S_2$. Constraints (A2)-(A10) corresponds to operation and security requirements of area $S_1$. Specifically, constraints (A2)-(A7) are local operation constraints for the generator in $S_1$ and constraints (A8)-(A9) are local network security constraints, both of which only include local variables in area S₁. On the other hand, constraint (A10) is a coupling constraint which includes tie-line power flow variable $PL_{23,t}$. Similarly, for area S₂, constraints (A11)-(A16) are local generator operation constraints, (A17) are local network security constraints, and (A18) are coupling constraints which include tie-line power flow variable $PL_{23,t}$.

As shown in Fig. 10(b), the original system can be decomposed into two independent sub-systems after duplicating bus 2/3 as bus 2'/3' in adjacent sub-system S₂/S₁. NCUC problems for the two sub-systems are presented as follows, where $\widetilde{PL}_{23',t}$ and $\widetilde{PL}_{2'3,t}$ are respectively power flow variables of the same tie-line at time $t$ in S₁ and S₂ (i.e., $\widetilde{PL}_{23',t} = \widetilde{PL}_{2'3,t}$).

| NCUC problem for S₁: | NCUC problem for S₂: |
|---|---|
| $\min \sum_{t=1}^{T}(F_1(P_{1,t}) + C_{nl1} \cdot I_{1,t} + SU_{1,t} + SD_{1,t})$ | $\min \sum_{t=1}^{T}(F_2(P_{2,t}) + C_{nl2} \cdot I_{2,t} + SU_{2,t} + SD_{2,t})$ |
| s.t. Constraints (A1)-(A12) | s.t. Constraints (A14)-(A23) |
| $PL_{12,t} - \widetilde{PL}_{23',t} = P_{d1,t}, \forall t \in T$ | $\widetilde{PL}_{2'3,t} - PL_{34,t} = P_{d2,t}, \forall t \in T$ |
| $\widetilde{PL}_{23',t} = \frac{(\theta_{2,t} - \theta_{3',t})}{X_{23}}, \forall t \in T$ | $\widetilde{PL}_{2'3,t} = \frac{(\theta_{2',t} - \theta_{3,t})}{X_{23}}, \forall t \in T$ |
| $PL_{23}^{min} \leq \widetilde{PL}_{23',t} \leq PL_{23}^{max}, \forall t \in T$ | $PL_{23}^{min} \leq \widetilde{PL}_{2'3,t} \leq PL_{23}^{max}, \forall t \in T$ |
| $\widetilde{PL}_{23',t} = Z_{1,t}, \forall t \in T$ | $\widetilde{PL}_{2'3,t} = Z_{1,t}, \forall t \in T$ |

Indexing $\widetilde{PL}_{23',t}$ as the second power flow variable in S₁ at time $t$ (i.e. $\widetilde{PL}_{S_1,t}(2) = \widetilde{PL}_{23',t}$) and $\widetilde{PL}_{2'3,t}$ as the second power flow variable in S₂ at time $t$ (i.e. $\widetilde{PL}_{S_2,t}(2) = \widetilde{PL}_{2'3,t}$), the mapping $G_t(n,w)=g$ between global variables $Z_{g,t}$ and duplicated local variables $PL_{n,t}(w)$ defined in Section II.B can be represented as $G_t(S_1,2)=1$ and $G_t(S_2,2)=1$ for each time $t \in T$. To match the consensus based ADMM form (7), coupling constraints $\widetilde{PL}_{23',t} = Z_{1,t}$ in NCUC problem 1 and $\widetilde{PL}_{2'3,t} = Z_{1,t}$ in NCUC problem 2 are relaxed into corresponded objective functions, as shown in the following sub-problems 1 and 2.

| Sub-problem 1: | Sub-problem 2: |
|---|---|
| $\min \sum_{t=1}^{T}\left[F_1(P_{1,t}) + C_{nl1} \cdot I_{1,t} + SU_{1,t} + SD_{1,t} + \lambda_{1t}(PL_{23',t} - Z_{1t}) + \rho(PL_{23',t} - Z_{1t})^2\right]$ | $\min \sum_{t=1}^{T}\left[F_2(P_{2,t}) + C_{nl2} \cdot I_{2,t} + SU_{2,t} + SD_{2,t} + \lambda_{2t}(PL_{2'3,t} - Z_{1t}) + \rho(PL_{2'3,t} - Z_{1t})^2\right]$ |
| s.t. Constraints (A1)-(A12) | s.t. Constraints (A14)-(A23) |
| $PL_{12,t} - \widetilde{PL}_{23',t} = P_{d1,t}, \forall t \in T$ | $\widetilde{PL}_{2'3,t} - PL_{34,t} = P_{d2,t}, \forall t \in T$ |
| $\widetilde{PL}_{23',t} = \frac{(\theta_{2,t} - \theta_{3',t})}{X_{23}}, \forall t \in T$ | $\widetilde{PL}_{2'3,t} = \frac{(\theta_{2',t} - \theta_{3,t})}{X_{23}}, \forall t \in T$ |
| $PL_{23}^{min} \leq \widetilde{PL}_{23',t} \leq PL_{23}^{max}, \forall t \in T$ | $PL_{23}^{min} \leq \widetilde{PL}_{2'3,t} \leq PL_{23}^{max}, \forall t \in T$ |

Finally, the consensus based ADMM procedure (i.e., Algorithm 1 in Section II.B) for this 2-area NCUC problem can be expressed as follows:

(1) Initializes $\lambda_{1,t}$, $\lambda_{2,t}$ and $Z_{1,t}$.
(2) S₁ solves sub-problem 1 and S₂ solves sub-problem 2 in a parallel manner.
(3) S₁ sends its $\widetilde{PL}_{23',t}$ solution to S₂, and S₂ sends its $\widetilde{PL}_{2'3,t}$ solution to S₁.
(4) Both areas S₁ and S₂ update $Z_{1,t}$ via equation $Z_{1,t} = \frac{(\widetilde{PL}_{23',t} + \widetilde{PL}_{2'3,t})}{2}, \forall t \in T$.
(5) S₁ updates $\lambda_{1,t}$ and S₂ updates $\lambda_{2,t}$ through the following two equations respectively,
$$\lambda_{1,t} = \lambda_{1,t} + \rho \cdot (\widetilde{PL}_{23',t} - Z_{1,t}), \forall t \in T$$
$$\lambda_{2,t} = \lambda_{2,t} + \rho \cdot (\widetilde{PL}_{2'3,t} - Z_{1,t}), \forall t \in T$$
(6) If the stopping criterion (7e) is satisfied, the algorithm terminates; otherwise, go back

to Step (2).